\documentclass[12pt,leqno]{article}
\usepackage{amssymb}
\usepackage[mathscr]{eucal}
\usepackage{amsmath,amssymb,latexsym,theorem,bbm} %
\usepackage{color}
\usepackage{appendix}

\newcommand{\NN}{\mathbb{N}}
\newcommand{\QQ}{\mathbb{Q}}
\newcommand{\RR}{\mathbb{R}}

\newcommand{\bone}{{\boldsymbol{1}}}

\newcommand{\cB}{{\mathcal B}}

\newcommand{\cP}{{\mathcal P}}

\newcommand{\dd}{\mathrm{d}}
\newcommand{\ee}{\mathrm{e}}

\newcommand{\EE}{\operatorname{\mathbb{E}}}
\newcommand{\PP}{\operatorname{\mathbb{P}}}

\newcommand{\tS}{\widetilde{S}}

\renewcommand{\leq}{\leqslant}
\renewcommand{\geq}{\geqslant}

\newcommand{\bbone}{\mathbbm{1}}

\newcommand{\proofend}{\hfill\mbox{$\Box$}}

\numberwithin{equation}{section}

\theoremstyle{change} \theorembodyfont{\em}
\newtheorem{Thm}{Theorem.}[section]

\newtheorem{Pro}[Thm]{Proposition.}

\newtheorem{Def}[Thm]{Definition.}

\theorembodyfont{\rm}
\newtheorem{Rem}[Thm]{Remark.}

\begin{document}

\begin{center}
 {\bfseries\Large A new example for a proper scoring rule}

\vspace*{3mm}

{\sc\large
  M\'aty\'as $\text{Barczy}^{*}$}

\end{center}

\vskip0.2cm

\noindent
 * MTA-SZTE Analysis and Stochastics Research Group,
   Bolyai Institute, University of Szeged,
   Aradi v\'ertan\'uk tere 1, H--6720 Szeged, Hungary.

\noindent e-mail: barczy@math.u-szeged.hu (M. Barczy).

\vskip0.2cm


\renewcommand{\thefootnote}{}
\footnote{\textit{2020 Mathematics Subject Classifications\/}:
  62C05, 62C99}
\footnote{\textit{Key words and phrases\/}:
  scoring rule, properization,  weighted Continuous Ranked Probability Scoring rule.}
\vspace*{0.2cm}
\footnote{M\'aty\'as Barczy is supported by the J\'anos Bolyai Research Scholarship
 of the Hungarian Academy of Sciences.}

\vspace*{-10mm}

\begin{abstract}
We give a new example for a proper scoring rule motivated by the form of Anderson--Darling distance of distribution functions
 and Example 5 in Brehmer and Gneiting (2020).
\end{abstract}

\section{Introduction and a new proper scoring rule}

Decisions based on accurate (probabilistic) forecasts of quantities of interest are very important in practice,
 since they affect our daily life very much.
For example, meteorology involves forecasting of temperature and wind speed, in hydrology it is important to predict water level,
 and in time series analysis, to forecast future values of some time series modelled, e.g., by some autoregressive moving average process.
Such a quantity in question is usually modelled by a random variable having an unknown distribution function, and, in general,
 several forecasting distribution functions are proposed by the practitioners, so it is a challenging and important task is
 to determine which one is the best and in which sense.
A so-called scoring rule assigns a score based on the forecasted distribution function and the realized observations, see, e.g., Gneiting and Raferty (2007)
 and David and Musio (2014).
More precisely, following the setup and notations of Brehmer and Gneiting (2020),
 let \ $\Omega$ \ be a non-empty set, \ $\cB$ \ be a \ $\sigma$-algebra on \ $\Omega$,
 and \ $\cP$ \ be a convex set of probability measures on \ $(\Omega,\cB)$.
\ A {\sl scoring rule} is an extended real valued function \ $S:\cP\times \Omega\to\RR\cup\{-\infty, \infty\}$ \ such that
 \[
   S(\PP,\QQ):= \int_\Omega S(\PP,\omega)\, \QQ(\dd\omega)
 \]
 is well-defined for all \ $\PP,\QQ\in\cP$ \ (in particular, for all \ $\PP\in\cP$, \ the mapping \ $\Omega\ni \omega\mapsto S(\PP,\omega)$ \ is measurable).
A scoring rule \ $S$ \ is called {\sl proper} relative to \ $\cP$, \ if
 \begin{align}\label{help6}
    S(\QQ,\QQ)\leq S(\PP,\QQ) \qquad \text{for all \ $\PP,\QQ\in \cP$.}
 \end{align}
A scoring rule \ $S$ \ is called {\sl strictly proper} relative to \ $\cP$, \ if it is proper relative to \ $\cP$, \
 and for any \ $\PP,\QQ\in\PP$, \ the equality \ $S(\QQ,\QQ) = S(\PP,\QQ)$ \ implies \ $\QQ=\PP$.
\ Note that in information theory a similar inequality to \eqref{help6} appears, namely, if \ $\xi$ \ and \ $\eta$ \ are discrete random variables
 having finite ranges, then \ $H(\xi,\xi)=H(\xi)\leq H(\xi,\eta)$, \ where \ $H(\xi,\xi)$ \ and \ $H(\xi,\eta)$ \ denote the Shannon entropy
 of \ $(\xi,\xi)$ \ and \ $(\xi,\eta)$, \ respectively.
In fact, the Shannon entropy \ $H(\xi) = -\sum_{i=1}^n p_i \log_2(p_i)$ \ of a discrete random variable \ $\xi$ \ having range
 \ $\{x_1,\ldots,x_n\}\subset\RR$ \ with some \ $n\in\NN$ \ and having distribution \ $\PP_\xi(\{x_i\}) = p_i$, \ $i=1,\ldots,n$, \
 coincides with \ $S_{\mathrm{ent}}(\PP_\xi,\PP_\xi) = \int_\RR S_{\mathrm{ent}}(\PP_\xi,\omega)\,\PP_\xi(\dd\omega)$, \ where
 \ $S_{\mathrm{ent}}(\PP_\xi,\omega):=-\log_2(p_i)$ \ if \ $\omega=x_i$ \ with some \ $i\in\{1,\ldots,n\}$ \ and
 \ $S_{\mathrm{ent}}(\PP_\xi,\omega):=0$ \ otherwise.
Here \ $S_{\mathrm{ent}}$ \  is nothing else but the so-called logarithmic score, for more details, see, e.g.,
 Gneiting and Raftery (2007, Example 3 and Section 4).

Next, we give an interpretation of the inequality \eqref{help6} in case of \ $\Omega$ \ is the real line \ $\RR$ \ and
 \ $\cB$ \ is the Borel \ $\sigma$-algebra \ $\cB(\RR)$ \ on \ $\RR$.
\ Suppose that we believe that a real-valued random variable \ $X$ \  has a distribution \ $\QQ$ \ on \ $(\RR,\cB(\RR))$,
\ and that the penalty for quoting some predictive distribution \ $\PP$ \ on \ $(\RR,\cB(\RR))$ \ for a realization
 \ $\omega\in\RR$ \ is \ $S(\PP,\omega)\in\RR\cup\{-\infty, \infty\}$.
\ So if our quoted distribution for \ $X$ \ is \ $\PP$, \ then the expected value of our penalty is
 \ $\EE(S(\PP,X)) = \int_\RR S(\PP,\omega)\, \QQ(\dd\omega) = S(\PP,\QQ)$.
\ Based on principles of decision theory, we should choose our quoting distribution \ $\PP$ \ in order to minimize
 the expected penalty \ $S(\PP,\QQ)$, \ and inequality \eqref{help6} says that \ $\QQ$ \ is such an optimal choice.
If \ $S$ \ is strictly proper, then \ $\QQ$ \ is the unique optimal choice.

It is known that a scoring rule satisfying some kind of regularity condition (see \eqref{help2}) can be properized in the sense that
 it can be modified in a way that it becomes a proper scoring rule, see, e.g., Theorem 1 in Brehmer and Gneiting (2020), which we recall below.

\begin{Thm}\label{Thm_properization}
Let \ $S:\cP\times \Omega\to\RR\cup\{-\infty, \infty\}$ \ be a scoring rule.
Suppose that for every \ $\PP\in\cP$ \ there exists a probability measure \ $\PP^*\in\cP$ \ such that
 \begin{align}\label{help2}
   S(\PP^*,\PP) \leq S(\QQ,\PP)\qquad \text{for all \ $\QQ\in \cP$.}
 \end{align}
Then the function \ $S^*:\cP\times\Omega\to \RR\cup\{-\infty, \infty\}$ \ defined by
 \[
   S^*(\PP,\omega):=S(\PP^*,\omega),\qquad \PP\in\cP, \; \omega\in\Omega,
 \]
 is a proper scoring rule.
\end{Thm}

In all what follows, let \ $\Omega:=\RR$ \ and \ $\cB$ \ be the Borel $\sigma$-algebra on \ $\RR$,
 \ and a probability measure \ $\PP\in\cP$ \ is identified with the function \ $\RR\ni x \mapsto \PP((-\infty,x))$,
 \ which is nothing else but the distribution function of the random variable \ $\Omega\ni\omega\mapsto \omega$ \
 with respect to the probability measure \ $\PP$.
\ Here we remind the readers that we use the previous definition of a distribution function instead of
 \ $\RR\ni x\mapsto \PP((-\infty,x])$ \ (which is also common in the literature).
In notation, instead of \ $\PP((-\infty, x))$ \ we will write \ $\PP(x)$, \ where \ $x\in\RR$.

A commonly used scoring rule is the so-called weighted Continuous Ranked Probability Scoring rule (wCRPS)
 defined by
 \[
     \mathrm{wCRPS}(\PP,y):=\int_{\RR} (\PP(x) - \bone_{\{y<x\}})^2 w(x)\,\dd x, \qquad \PP\in\cP, \; y\in\RR,
 \]
 where \ $w:\RR\to(0,\infty)$ \ is a given measurable function (called weight function as well).
In the special case \ $w(x)=1$, \ $x\in\RR$, \ wCRPS is nothing else but the Continuous Ranked Probability Scoring rule (CRPS),
 see, e.g., Gneiting and Raferty (2007, Section 4.2).
Sometimes, wCRPS and CRPS is simply called weighted Continuous Ranked Probability Score
 and Continuous Ranked Probability Score, respectively.
By choosing weight functions in an appropriate way the center of (one of the) tails of the range of distribution functions can be emphasized.
For more details on the role of weight functions and examples for some commonly used weight functions, see Gneiting and Ranjan (2011, page 415 and Table 4).
These scoring rules are commonly used in practice, see, e.g., the very recent work of Baran, Hemri and El Ayari (2019).

Recently, for any \ $\alpha>0$, \ Brehmer and Gneiting (2020, Example 5) have introduced a scoring rule
 \ $S_\alpha: \cP\times \RR\to [0,\infty]$ \ given by
 \begin{align}\label{help3}
  S_\alpha(\PP,y):= \int_{\RR} \vert \PP(x) - \bone_{\{y<x\}}\vert^\alpha\,\dd x, \qquad \PP\in\cP, \; y\in\RR.
 \end{align}
For \ $\alpha=2$, \ it gives back the scoring rule CRPS.
Using Theorem \ref{Thm_properization}, Brehmer and Gneiting (2020) have shown that
 the function \ $S_\alpha^*:\cP\times \RR\to [0,\infty]$,
 \[
    S_\alpha^*(\PP,y):=S_\alpha(\PP^*,y), \qquad \PP\in\cP, \; y\in\RR,
 \]
 is a proper scoring rule, where the mapping \ $\cP\ni \PP\mapsto \PP^*\in\cP$ \ is given by
 \[
   \PP^*(x):=\left(1+\left(\frac{1-\PP(x)}{\PP(x)}\right)^{\frac{1}{\alpha-1}} \right)^{-1}\bone_{\{\PP(x)>0\}}, \qquad \PP\in\cP,
   \qquad \text{in case of \ $\alpha>1$,}
 \]
 and \ $\PP^*$ \ is the distribution function of the Dirac measure concentrated at a median of \ $\PP$ \ in case of \ $\alpha\in(0,1]$.
\ If \ $\alpha=1$ \ and there is more than one median of \ $\PP$, \ then there are other choices for \ $\PP^*$.
\ In case of \ $\alpha>1$, \ the mapping \ $\cP\ni \PP\mapsto \PP^*\in\cP$ \ is injective.
\ Further, in some cases the proper scoring rule \ $S_\alpha$ \ can be made a strictly proper scoring rule.
For example, if \ $\alpha\in(1,2]$, \  then \ $S_\alpha$ \ restricted to \ $\cP_1\times \RR$ \ is a strictly proper
 scoring rule relative to \ $\cP_1$, \ where \ $\cP_1$ \ denotes the set of probability measures on \ $(\RR,\cB(\RR))$ \ with finite first
 moment, see Brehmer and Gneiting (2020, Example 5).

Motivated by \eqref{help3} and the form of Anderson--Darling distance of distribution functions (see, e.g., Anderson and Darling (1954)
 or Deza and Deza (2013, page 237), we introduce a new scoring rule.
Let \ $\cP_{(0,1)}$ \ be the set of distribution functions taking values in \ $(0,1)$.
Then \ $\cP_{(0,1)}$ \ is a convex subset of \ $\cP$.

\begin{Def}
For \ $\alpha>0$ \ and a measurable function \ $w:\RR\to (0,\infty)$,
 \ let \ $\tS_{\alpha,w}:\cP_{(0,1)}\times \RR\to [0,\infty] $,
 \[
  \tS_{\alpha,w}(\PP,y):=\int_\RR \frac{\vert \PP(x) - \bbone_{\{y< x\}}\vert^{2\alpha}}{\PP(x)^\alpha(1-\PP(x))^\alpha}w(x)\,\dd x,
  \qquad \PP\in\cP_{(0,1)}, \; y\in\RR.
 \]
\end{Def}

\begin{Pro}\label{Pro_new_scoring_rule}
For each \ $\alpha>0$ \ and for each measurable function \ $w:\RR\to (0,\infty)$, \
 the function \ $\tS_{\alpha,w}^*:\cP_{(0,1)}\times \RR\to [0,\infty]$,
 \[
 \tS_{\alpha,w}^*(\PP,y):=\tS_{\alpha,w}({\widetilde\PP}^*,y), \qquad \PP\in\cP_{(0,1)}, \;\; y\in\RR,
 \]
 is a proper scoring rule, where the mapping \ $\cP_{(0,1)}\ni \PP\mapsto \widetilde{\PP}^*\in\cP_{(0,1)}$ \ is given by
 \begin{align}\label{help1}
  \widetilde{\PP}^*(x) := \left(1 + \left(\frac{1-\PP(x)}{\PP(x)} \right)^{\frac{1}{2\alpha}} \right)^{-1},\qquad x\in\RR.
 \end{align}
Further, for any \ $\PP\in\cP_{(0,1)}$ \ and \ $y\in\RR$,
 \begin{align}\label{help5}
   \tS_{\alpha,w}({\widetilde\PP}^*,y)
      = \int_\RR \frac{ \vert \bbone_{(y,\infty)}(x) - \PP(x) \vert^{\frac{1}{2}} }{ \vert 1 - \bbone_{(y,\infty)}(x) - \PP(x) \vert^{\frac{1}{2}}}w(x)\,\dd x,
 \end{align}
 and
 \begin{align}\label{help7}
   \tS_{\alpha,w}({\widetilde\PP}^*,\PP)
       = \int_{\RR} \tS_{\alpha,w}({\widetilde\PP}^*,y)\, \PP(\dd y)
      = 2\int_\RR \big(\PP(x) (1-\PP(x))\big)^{\frac{1}{2}}w(x)\,\dd x.
 \end{align}
\end{Pro}

The proof of Proposition \ref{Pro_new_scoring_rule} can be found in Section \ref{Sec_Proof}.
In the next remark we give an example, where a restriction of \ $\tS_{\alpha,w}^*$ \ leads to a {\sl strictly} proper scoring rule.

\begin{Rem}\label{Rem1}
\noindent{(i)}
Let \ $\alpha>0$ \ and \ $w:\RR\to (0,\infty)$, \ $w(x):=1$, \ $x\in\RR$.
\ Let \ ${\widehat \cP}_{(0,1)}$ \ be a subclass of \ $\cP_{(0,1)}$ \ satisfying the following two properties:
 \ $\PP$ \  is continuous for any \ $\PP\in {\widehat \cP}_{(0,1)}$, \ and \ $\tS_{\alpha,w}(\widetilde\PP^*,\PP)$ \ is finite for
 any \ $\PP\in {\widehat \cP}_{(0,1)}$.
\ Then \ $\tS_{\alpha,w}^*$ \ restricted to \ ${\widehat \cP}_{(0,1)}\times\RR$ \ is strictly proper relative to
 \ ${\widehat \cP}_{(0,1)}$, \ i.e., it is proper relative to \ ${\widehat \cP}_{(0,1)}$, \
 and for any \ $\PP,\QQ\in {\widehat \cP}_{(0,1)}$, \ the equality
 \ $\tS_{\alpha,w}^*(\QQ,\QQ) = \tS_{\alpha,w}^*(\PP,\QQ)$ \ implies \ $\QQ=\PP$.
\ For a proof, see  Section \ref{Sec_Proof}.

{(ii)}
If \ $w(x)=1$, \ $x\in\RR$, \ and \ $\PP(x):= \ee^{-\ee^{-x}}$, \ $x\in\RR$ \ (Gumbel distribution),
 or
 \[
    \PP(x):=\begin{cases}
             \frac{1}{2}\ee^{x} & \text{if \ $x\leq0$,}\\
             1-\frac{1}{2}\ee^{-x} & \text{if \ $x> 0$,}
            \end{cases}
 \]
 (Laplace distribution), then \ $\tS_{\alpha,w}({\widetilde\PP}^*,\PP)$ \ is finite.
\proofend
\end{Rem}

In the next remark we initiate two other scoring rules.

\begin{Rem}
One may try to investigate the properties of the scoring rules
 \[
   \int_\RR \vert \PP(x) - \bbone_{\{y< x\}}\vert^{2\alpha}\,\PP(\dd x), \qquad y\in\RR, \;\; \PP\in\cP,
 \]
 and
 \[
   \int_\RR \frac{\vert \PP(x) - \bbone_{\{y< x\}}\vert^{2\alpha}}{\PP(x)^\alpha(1-\PP(x))^\alpha}\,\PP(\dd x), \qquad y\in\RR, \;\; \PP\in\cP_{(0,1)},
 \]
 where \ $\alpha>0$.
\ The second one with \ $\alpha=1$ \ is nothing else but the Anderson-Darling distance of the distribution functions \ $\PP(x)$, $x\in\RR$, \
 and \ $\bbone_{\{y< x\}}$, \ $x\in\RR$.
\ For these scoring rules we were not able to derive similar results as in Proposition \ref{Pro_new_scoring_rule}.
\proofend
\end{Rem}

\section{Proofs}\label{Sec_Proof}

\noindent{\bf Proof of Proposition \ref{Pro_new_scoring_rule}.}
The technique of our proof is similar to that of Example 5 in Brehmer and Gneiting (2020), namely,
 we will use Theorem \ref{Thm_properization}.
\ Fix \ $\alpha>0$, \ $\PP\in\cP_{(0,1)}$ \ and a measurable function \ $w:\RR\to(0,\infty)$.
\ Then for all \ $\QQ\in\cP_{(0,1)}$, \ by Tonelli's theorem,
 \begin{align}\label{help8}
  \begin{split}
   &\tS_{\alpha,w}(\QQ,\PP)
             = \int_\RR \tS_{\alpha,w}(\QQ,y)\, \PP(\dd y)
             = \int_\RR\left( \int_\RR  \frac{\vert \QQ(x) - \bbone_{\{ y< x\}}\vert^{2\alpha}}{\QQ(x)^\alpha(1-\QQ(x))^\alpha}\, w(x)\,\dd x\right)\PP(\dd y)\\
           & = \int_\RR\left( \int_\RR \frac{\vert\QQ(x) - \bbone_{\{ y< x\}}\vert^{2\alpha}}{\QQ(x)^\alpha(1-\QQ(x))^\alpha} \,\PP(\dd y) \right)w(x) \dd x \\
           & = \int_\RR\left( \int_{\{y< x\}} \frac{\vert\QQ(x) - 1\vert^{2\alpha}}{\QQ(x)^\alpha(1-\QQ(x))^\alpha} \,\PP(\dd y) \right) w(x)\dd x\\
           &\phantom{=\;} + \int_\RR\left( \int_{\{y\geq x\}} \frac{\QQ(x)^{2\alpha}}{\QQ(x)^\alpha(1-\QQ(x))^\alpha} \,\PP(\dd y) \right) w(x)\dd x \\
           & =  \int_\RR\left( \left(\frac{1-\QQ(x)}{\QQ(x)}\right)^\alpha \PP(x) + \left( \frac{\QQ(x)}{1-\QQ(x)}\right)^\alpha (1-\PP(x)) \right) w(x)\dd x.
  \end{split}
 \end{align}
For fixed \ $x\in\RR$, \ let us introduce the function \ $g_{x,\PP}:(0,1)\to\RR$,
 \[
   g_{x,\PP}(q):= \left[\left(\frac{1-q}{q}\right)^\alpha \PP(x) + \left(\frac{q}{1-q}\right)^\alpha(1-\PP(x))\right]w(x), \qquad q\in(0,1).
 \]
One can calculate that for any \ $q\in(0,1)$,
 \begin{align*}
   &g_{x,\PP}'(q)  =  \alpha \left(\frac{1-q}{q}\right)^{\alpha-1} \left(-\frac{1}{q^2}\right) \PP(x) w(x)
                    + \alpha \left(\frac{q}{1-q}\right)^{\alpha-1} \frac{1}{(1-q)^2}(1-\PP(x))w(x),\\[1mm]
   &g_{x,\PP}''(q) =  \alpha \left(\frac{1-q}{q}\right)^{\alpha-2}  \frac{\alpha+1-2q}{q^4} \PP(x)w(x)
                      + \alpha \left(\frac{q}{1-q}\right)^{\alpha-2}  \frac{\alpha-1+2q}{(1-q)^4} (1-\PP(x))w(x).
 \end{align*}

If \ $\alpha\geq 1$, \ then \ $g_{x,\PP}''(q)>0$, \ $q\in(0,1)$, \
 and hence the function \ $g_{x,\PP}$ \ is strictly convex on \ $(0,1)$, \
 and its unique minimum is attained at \ $q_{x,\PP}^*\in(0,1)$, \ which satisfies the equation \ $g_{x,\PP}'(q_{x,\PP}^*)=0$. \
One can calculate that
 \begin{align}\label{help4}
   q_{x,\PP}^* = \left(1 + \left( \frac{1-\PP(x)}{\PP(x)} \right)^{\frac{1}{2\alpha}} \right)^{-1}=: {\widetilde\PP}^*(x), \qquad x\in\RR.
 \end{align}

Next, we show that for all \ $\alpha>0$, \ the function \ $g_{x,\PP}$ \ attains its minimum at \ $g_{x,\PP}^*$ \ given in \eqref{help4}.
Since \ $g_{x,\PP}'(q_{x,\PP}^*) = 0$, \ for this, it is enough to check that \ $g_{x,\PP}''(q_{x,\PP}^*) >0$.
\ Since for all \ $q\in(0,1)$,
 \[
   g_{x,\PP}''(q) = \alpha w(x)\left(\frac{1-q}{q}\right)^{\alpha-2}
                    \left[
                         \frac{\alpha+1-2q}{q^4} \PP(x) + \left(\frac{q}{1-q}\right)^{2(\alpha-2)} \frac{\alpha-1+2q}{(1-q)^4} (1-\PP(x))
                    \right]  ,
 \]
 we have
 \begin{align*}
      &g_{x,\PP}''(q_{x,\PP}^*) \\
      & = \alpha w(x) \left( \frac{1-\PP(x)}{\PP(x)} \right)^{\frac{\alpha-2}{2\alpha}}
          \left[
                \frac{\alpha+1-2q_{x,\PP}^*}{(q_{x,\PP}^*)^4} \PP(x) +  \left( \frac{\PP(x)}{1-\PP(x)} \right)^{1-\frac{2}{\alpha}} \frac{\alpha-1+2q_{x,\PP}^*}{(1-q_{x,\PP}^*)^4} (1-\PP(x))
          \right] \\
      & = \alpha w(x) \PP(x) \left( \frac{1-\PP(x)}{\PP(x)} \right)^{\frac{\alpha-2}{2\alpha}}
          \left[
                \frac{\alpha+1-2q_{x,\PP}^*}{(q_{x,\PP}^*)^4}  +  \left( \frac{1-\PP(x)}{\PP(x)} \right)^{\frac{4}{2\alpha}} \frac{\alpha-1+2q_{x,\PP}^*}{(1-q_{x,\PP}^*)^4}
          \right]\\
      & = \alpha w(x) \PP(x) \left( \frac{1-\PP(x)}{\PP(x)} \right)^{\frac{\alpha-2}{2\alpha}}
           \left[
                \frac{\alpha+1-2q_{x,\PP}^*}{(q_{x,\PP}^*)^4}  +  \left( \frac{1- q_{x,\PP}^*}{q_{x,\PP}^*}  \right)^4 \frac{\alpha-1+2q_{x,\PP}^*}{(1-q_{x,\PP}^*)^4}
          \right]\\
      & = \alpha w(w)\PP(x) \left(\frac{1}{q_{x,\PP}^*}-1\right)^{\alpha-2} \frac{2\alpha}{(q_{x,\PP}^*)^4}
      >0,
 \end{align*}
 as desired.
One can easily check that \ ${\widetilde\PP}^*\in\cP_{(0,1)}$, \ i.e., \ ${\widetilde\PP}^*$ \ is a distribution function with values in \ $(0,1)$.
\ Note also that \eqref{help1} shows that the mapping \ $\cP_{(0,1)}\ni \PP\mapsto \widetilde{\PP}^*$ \ is injective.
All in all, condition \eqref{help2} of Theorem \ref{Thm_properization} with \ $\cP:=\cP_{(0,1)}$ \ is satisfied, i.e.,
 for every \ $\PP\in\cP_{(0,1)}$ \ there exists \ $\PP^*\in\cP_{(0,1)}$ \ such that
 \[
  \tS_{\alpha,w}({\widetilde\PP}^*,\PP) \leq \tS_{\alpha,w}(\QQ,\PP), \qquad \QQ\in\cP_{(0,1)},
 \]
 and, by Theorem \ref{Thm_properization}, we have the first part of the assertion.
Further, for any \ $\PP\in\cP_{(0,1)}$ \ and \ $y\in\RR$,
 \begin{align*}
  \tS_{\alpha,w}({\widetilde\PP}^*,y)
    & = \int_{\{y<x\}} \frac{\vert {\widetilde\PP}^*(x) - 1\vert^{2\alpha}}{{\widetilde\PP}^*(x)^\alpha(1-{\widetilde\PP}^*(x))^\alpha} w(x)\,\dd x
       + \int_{\{y\geq x\}} \frac{{\widetilde\PP}^*(x)^{2\alpha}}{{\widetilde\PP}^*(x)^\alpha(1-{\widetilde\PP}^*(x))^\alpha} w(x)\,\dd x \\
    & = \int_{\{y<x\}} \left( \frac{1}{{\widetilde\PP}^*(x)} - 1 \right)^\alpha w(x)\,\dd x
        + \int_{\{y\geq x\}} \left( \frac{1}{{\widetilde\PP}^*(x)} - 1 \right)^{-\alpha}  w(x)\,\dd x\\
    & = \int_{\{y<x\}} \left( \frac{1-\PP(x)}{\PP(x)} \right)^{\frac{1}{2}} w(x)\,\dd x
        + \int_{\{y\geq x\}} \left( \frac{\PP(x)}{1-\PP(x)}  \right)^{\frac{1}{2}}  w(x)\,\dd x\\
    & = \int_\RR \frac{ \vert \bbone_{(y,\infty)}(x) - \PP(x) \vert^{\frac{1}{2}} }{ \vert 1 - \bbone_{(y,\infty)}(x) - \PP(x) \vert^{\frac{1}{2}}} w(x)\,\dd x,
 \end{align*}
 as desired.

Finally, by \eqref{help8} and the definition of the function \ $g_{x,\PP}$, \ we have
 \begin{align*}
   \tS_{\alpha,w}({\widetilde\PP}^*,\PP)  = \int_\RR g_{x,\PP}({\widetilde\PP}^*(x))\,\dd x = 2\int_\RR \big(\PP(x) (1-\PP(x))\big)^{\frac{1}{2}}w(x)\,\dd x ,
 \end{align*}
 since
 \begin{align*}
  g_{x,\PP}({\widetilde\PP}^*(x))
    &= \left(\frac{1}{{\widetilde\PP}^*(x)} -1\right)^\alpha\PP(x)w(x) + \left(\frac{1}{{\widetilde\PP}^*(x)} -1\right)^{-\alpha}(1-\PP(x))w(x)\\
    &= \left(\frac{1-\PP(x)}{\PP(x)}\right)^{\frac{1}{2}}\PP(x)w(x) +  \left(\frac{1-\PP(x)}{\PP(x)}\right)^{-\frac{1}{2}} (1-\PP(x))w(x)\\
    &= 2(\PP(x)(1-\PP(x)))^{\frac{1}{2}}w(x), \qquad x\in\RR.
 \end{align*}
\proofend

\smallskip

\noindent{\bf Proof of the example given in part (i) of Remark \ref{Rem1}.}
By Proposition \ref{Pro_new_scoring_rule}, \ $\tS_{\alpha,w}^*:\cP_{(0,1)}\times \RR \to [0,\infty]$ \ is a proper scoring rule relative to
 \ $\cP_{(0,1)}$, \ so, especially, \ $\tS_{\alpha,w}^*(\QQ,\QQ) \leq \tS_{\alpha,w}^*(\PP,\QQ)$ \ for all \ $\PP,\QQ\in\widehat\cP_{(0,1)}$ \
 yielding that the restriction \ $\tS_{\alpha,w}^{*,r}: \widehat\cP_{(0,1)} \times \RR \to [0,\infty]$ \ is a proper scoring rule relative to
 \ $\widehat\cP_{(0,1)}$.
\ It remains to check the strict property of \ $\tS_{\alpha,w}^{*,r}: \widehat\cP_{(0,1)} \times \RR \to[0,\infty]$.
\ Let \ $\PP,\QQ\in \widehat\cP_{(0,1)}$ \ be such that \ $\tS_{\alpha,w}^{*,r}(\QQ,\QQ) = \tS_{\alpha,w}^{*,r}(\PP,\QQ)$.

Using \eqref{help8} and \eqref{help1} we have
 \begin{align*}
   \tS_{\alpha,w}^{*,r}(\PP,\QQ)
            & = \int_\RR \tS_{\alpha,w}^{*,r}(\PP,y)\, \QQ(\dd y)
              = \int_\RR \tS_{\alpha,w}^{*}(\PP,y)\, \QQ(\dd y)
              = \int_\RR \tS_{\alpha,w}({\widetilde\PP}^*,y)\, \QQ(\dd y)
              = \tS_{\alpha,w}({\widetilde\PP}^*,\QQ) \\
            & = \int_\RR\left( \left(\frac{1-{\widetilde\PP}^*(x)}{{\widetilde\PP}^*(x)}\right)^\alpha \QQ(x)
                 + \left( \frac{{\widetilde\PP}^*(x)}{1-{\widetilde\PP}^*(x)}\right)^\alpha (1-\QQ(x)) \right) \dd x\\
            & = \int_\RR\left( \left(\frac{1-\PP(x)}{\PP(x)}\right)^{\frac{1}{2}} \QQ(x)
                 + \left( \frac{1-\PP(x)}{\PP(x)}\right)^{-\frac{1}{2}} (1-\QQ(x)) \right) \dd x,
 \end{align*}
 and similarly (or referring to \eqref{help7})
 \[
   \tS_{\alpha,w}^{*,r}(\QQ,\QQ) = \tS_{\alpha,w}({\widetilde\QQ}^*,\QQ)
                                 = 2 \int_\RR (\QQ(x)(1-\QQ(x)))^{\frac{1}{2}}\,\dd x.
 \]
Note that, by the inequality between the arithmetic mean and the geometric mean, for any \ $x\in\RR$ \ we have
 \begin{align}\label{help9}
   \big(\QQ(x)(1-\QQ(x))\big)^{\frac{1}{2}}
      \leq \frac{1}{2}
             \left(\left(\frac{1-\PP(x)}{\PP(x)}\right)^{\frac{1}{2}} \QQ(x)
                    + \left( \frac{1-\PP(x)}{\PP(x)}\right)^{-\frac{1}{2}} (1-\QQ(x)) \right),
 \end{align}
 and equality holds if and only if
 \begin{align}\label{help10}
   \left(\frac{1-\PP(x)}{\PP(x)}\right)^{\frac{1}{2}} \QQ(x)
      = \left( \frac{1-\PP(x)}{\PP(x)}\right)^{-\frac{1}{2}} (1-\QQ(x))
      \qquad \Longleftrightarrow \qquad \PP(x) = \QQ(x).
 \end{align}
The inequality \eqref{help9} directly shows that \ $\tS_{\alpha,w}^{*,r}(\QQ,\QQ) \leq \tS_{\alpha,w}^{*,r}(\PP,\QQ)$, \ $\PP,\QQ\in \widehat\cP_{(0,1)}$
 \ (which we already know, since \ $\tS_{\alpha,w}^{*,r}$ \ is a proper scoring rule relative to \ $\widehat\cP_{(0,1)}$), \ and using the
 inequality \eqref{help10} and the continuity of \ $\PP$ \ and \ $\QQ$, \ a standard measure theoretical argument
 yields that \ $\tS_{\alpha,w}^{*,r}(\QQ,\QQ) = \tS_{\alpha,w}^{*,r}(\PP,\QQ)$ \ holds if and only if \ $\QQ=\PP$.
\proofend

\section*{Acknowledgements}
I am grateful to Jonas Brehmer for his comments to preliminary versions of the paper that helped me a lot.
I thank S\'andor Baran for calling my attention to the paper of Gneiting and Ranjan (2011).
I would like to thank the referee for the comments that helped me to improve the paper.

\end{document}